%% file: CCirred.tex
\documentclass{amsart}

\input{preamble.tex}

\input{definitions.tex}
\begin{document}

\title[IC irreducibility]{Characteristic cycles of IH sheaves of simply laced minuscule Schubert varieties are irreducible} 
\author{Leonardo C.  Mihalcea}
\address{Department of Mathematics, Virginia Tech, 460 McBryde Hall, 225 Stanger St., Blacksburg VA 24061}
\email{lmihalce@math.vt.edu} 

\author{Rahul Singh}
\email{rahul.sharpeye@gmail.com}

\thanks{L.~C.~Mihalcea was supported in part by the NSF grant DMS-2152294 and the Simons Collaboration Grant 581675.}

\date{August 11, 2023} 
\maketitle

\begin{abstract} 
Let $G/P$ be a complex cominuscule flag manifold of type $E_6,E_7$.
We prove that each characteristic cycle of the intersection homology (IH) 
complex of a Schubert variety in $G/P$ is irreducible. The proof utilizes 
an earlier algorithm by the same authors which calculates local Euler obstructions, 
then proceeds by direct computer 
calculation using \texttt{Sage}.
This completes to the exceptional Lie types the characterization of 
irreducibility of IH sheaves of Schubert varieties in cominuscule $G/P$ 
obtained by Boe and Fu. As a by-product, we also obtain that the Mather classes, and the 
Chern-Schwartz-MacPherson classes of Schubert cells in cominuscule $G/P$ of type $E_6,E_7$, 
are strongly positive.
\end{abstract}

\section{Introduction}\label{sec:intro} Let $G$ be a complex semisimple Lie 
group with Weyl group $W$ and let $B \subset P \subset G$ be a Borel group, 
included in a parabolic subgroup. Denote by $X=G/P$ the 
associated generalized flag manifold. The Borel group acts naturally by left multiplication,
with finitely many $B$-orbits $X_w^\circ = Bw P/P$. Here $w$ varies in the set of
minimal length representatives of $W/W_P$, the quotient of $W$ by the Weyl group of $P$. The variety 
$X_w^\circ$ is called a Schubert cell, and its closure $X_w = \overline{X_w^\circ} \subset X$ a Schubert variety.    

If $Y \subset X$ is a closed irreducible algebraic
subvariety, the intersection homology (IH) characteristic cycle of $Y$ is a conic Lagrangian cycle
in the cotangent bundle of $X$. For the purposes of this paper, this is an element  
$IH(Y) \in H^{2 \dim X}_{\C^*} (T^*(X)) $, 
the middle $\C^*$-equivariant cohomology of $T^*(X)$, where $\C^*$ acts by fibrewise dilation.
The group $H^{2 \dim X}_{\C^*} (T^*(X))$ is a free abelian group with a basis given by the
fundamental classes $[T^*_{X_v}X]$ of the conormal spaces of Schubert varieties; see e.g. \cite{HTT}.
The conormal space is defined as the closure $T^*_{X_v}X = \overline{T^*_{X_v^\circ}X} \subset T^*X$. 
Consider the expansion   
\[ IH(Y) = \sum_v m_Y^v [T^*_{X_v} X] \quad \in H^{2 \dim X}_{\C^*} (T^*(X)) \/. \]
Here $m_Y^v$ are non-negative integers called the IH multiplicities of $Y$. 

A fundamental problem in geometric representation theory is to calculate the multiplicities
$m_{w,v}:=m_{X_w}^v$, in the special case when $Y= X_w$ is a Schubert variety \cite{KL:topological,kashiwara.tanisaki:characteristic,MR1084458,boe.fu,
evens.mirkovic:characteristic,braden:irred,williamson:reducible}. 
Besides their intrinsic interest, these multiplicities are 
related to problems in representation theory in characteristic $p$; see e.g. 
\cite{vilonen.williamson,saito:char}. 
A case of particular interest
is when $IH(X_w)$ is irreducible, or, equivalently, when $IH(X_w) = [T^*_{X_w} X]$. The following is the main result of this note, and it was announced in \cite{mihalcea.singh:conormal}.
\begin{theorem}\label{thm:intro} Let $G/P$ be a minuscule homogeneous space of Lie type $E_6$ or $E_7$.
Then for any Schubert variety $X_w \subset G/P$, the characteristic cycle $IH(X_w)$ is irreducible.\end{theorem}

Minuscule homogeneous spaces belong to the larger family of {\em cominuscule} spaces, and 
they have many common geometric and combinatorial properties.
We present a complete list of these below, and we refer to, 
e.g.,~ \cite{BCMP:qkchev} for more about their properties.
\begin{itemize}[leftmargin=*]
\item
the Grassmannian $\Gr(k,n)$ (type A);
\item
the Lagrangian Grassmannian $\mathrm{LG}(n,2n)$ (type C);
\item
the connected components of the orthogonal Grassmannian $\mathrm{OG}(n,2n)$ (type D);
\item
Quadrics: odd dimensional in type B, and even dimensional in type D;
\item
the Cayley plane (type $E_6$), and the Freudenthal variety (type $E_7$).
\end{itemize}
{
Further, the orthogonal Grassmannian $\mathrm{OG}(n-1,2n-1)$ (a cominuscule space of type B) is isomorphic to a connected component of $\mathrm{OG}(n,2n)$,}
hence can be identified with a cominuscule space of type D, 
see~\cite[p.~68]{fulton.pragacz:schubert} or e.g.~\cite[\S 3.4]{IMN:factorial}.
In this note we will focus on the Cayley plane and the Fruedenthal variety. 
The minuscule spaces are those from the simply laced types (A,D,E).
The last two spaces are those from \Cref{thm:intro}. 
\input{cotcomin_table-1.tex}

If $G/P$ is a Grassmann manifold, Bressler, Finkelberg and Lunts \cite{MR1084458}
proved that the characteristic cycles of the IH sheaves of Schubert varieties are irreducible.
This was further studied by Boe and Fu \cite{boe.fu},
who utilized methods from geometric analysis 
to recover the results from \cite{MR1084458}, and in addition
calculated the multiplicities $m_{w,v}$ of the cycles $IH(X_w)$
for the Schubert varieties $X_w$ in the maximal orthogonal Grassmannians, and for the quadrics.
Combining Boe and Fu's results together with \Cref{thm:intro} above gives the following rather general statement.
\begin{corollary}\label{thm:intro2} Let $G/P$ be any cominuscule space. Then the IH sheaf of each Schubert variety 
in $G/P$ is irreducible if and only if $G$ is simply laced.\end{corollary}
It is not difficult to exhibit reducible IC characteristic cycles in the non simply laced cases; examples go back to
\cite{kashiwara.tanisaki:characteristic}, see also \cite{boe.fu} and more recently the authors' paper 
\cite{mihalcea.singh:conormal}. 
If one leaves the family of cominuscule spaces,
the IH sheaves of Schubert varieties are reducible in general, even in simply laced cases. For instance, 
Kashiwara and Saito \cite{kashiwara.saito:geometric} show that all Schubert varieties in 
$\Fl(n)= \SL_n/B$ for $n \le 7$ have irreducible IH characteristic cycles; they also 
provide an example of a Schubert variety in $\Fl(8)$ where this characteristic cycle is reducible. 

At the heart of our proof are the type uniform formulae obtained by the authors in \cite{mihalcea.singh:conormal} for calculating Mather classes of Schubert varieties
in cominuscule spaces $G/P$. Utilizing this, and as observed e.g.~in \cite[\S 8.3, Eq.~(40)]{AMSS:shadows} the irreducibility
of $IH(X_w)$ is equivalent to the equality
\begin{equation}\label{E:main} e_{w,v} = P_{w,v}(1) \quad \forall v \le w \/.\end{equation}
Here on the left hand side is the {\em local Euler obstruction} of $X_w$ at a (torus) fixed point $v$, a subtle invariant of singularities defined by MacPherson \cite{macpherson:chern}; our results from \cite{mihalcea.singh:conormal} give (rather complicated) cohomological formulae for this invariant in the cominuscule cases, but in general its calculation seems to be wide open.
On the right hand side is the (parabolic) {\em Kazhdan-Lusztig (KL) polynomial} $P_{w,v}(q)$ evaluated at $q=1$; see, e.g., \cite{deodhar:geometricII}. Our proof is by direct calculation. We utilize the aforementioned formulae from \cite{mihalcea.singh:conormal}, and known algorithms to calculate KL polynomials in this situation, then we wrote a \textsc{Sage} code which
checks \Cref{E:main}. The main steps of the code are explained below, and 
the code is available at \cite{MScode}. A natural and interesting project would be to obtain a type independent mathematical proof of \Cref{E:main} in the minuscule cases.

The present calculations in Lie type $E_6, E_7$ complete the proofs of several conjectures about Euler obstructions, Mather classes and Chern-Schwartz-MacPherson classes. Although our contribution in this note is solely restricted to these types, for the convenience of the reader we state the consequences here in all Lie types. For that, we will briefly recall few basic facts about CSM and Mather classes below. More details, including precise definitions of the objects involved, can be found in e.g. \cite{mihalcea.singh:conormal}.  

\subsection{Outline of the algorithm}
Let $X$ be a complex algebraic variety. The 
{\em Chern-Schwartz-MacPherson} (CSM) class of a constructible
subset $\Omega \subset X$ is an element  
$\csm(\Omega) \in H_*(X)$ in the homology group of $X$.
This was defined by 
MacPherson \cite{macpherson:chern} in relation to a conjecture by Grothendieck 
and Deligne about Chern classes of singular spaces. 
The CSM classes satisfy 
$\csm(\Omega_1 \cup \Omega_2) = \csm(\Omega_1) + \csm(\Omega_2)$ 
for disjoint constructible sets $\Omega_1, \Omega_2$ and are functorial 
with respect to proper push-forwards. Furthermore, the CSM classes are determined 
by these properties together with the normalization property that $\csm(X) = c(T(X)) \cap [X]$
if $X$ is smooth.

In 
\cite{aluffi.mihalcea:csm,aluffi.mihalcea:eqcsm} Aluffi and the first named 
author calculated CSM classes of Schubert cells in any (generalized) 
flag manifold $G/P$. We found explicit algorithms 
to calculate the expansions of CSM classes into Schubert
classes in terms of the Demazure-Lusztig (DL) operators:
\begin{equation}\label{E:CSM-Schub} \csm(X_w^\circ) = \sum_{v \le w} a_{w,v} [X_v] \quad \in H_*(G/P) \/. \end{equation}
The coefficients $a_{w,v}$ are non-negative integers; this was proved by J. Huh
\cite{huh:csm} for Grassmannians (see also \cite{jones:csm,MR2949825,mihalcea:binomial})
and for arbitrary $G/P$ by Aluffi, Sch{\"u}rmann, Su and the first named author in \cite{AMSS:shadows}.
 
A related class, but much harder to calculate, is the {\em Mather class}
of a Schubert variety $\cMa(X_w) \in H_*(G/P)$. If $X_w$ is smooth, both 
$\csm(X_w)$ and $\cMa(X_w)$ are equal to $i_*(c(T X_w))$, where 
$i: X_w \hookrightarrow G/P$ is the inclusion. Historically, the Mather class 
was defined first, and it was utilized by MacPherson to define CSM classes.
From MacPherson's definition, the CSM classes and the Mather classes
are related by the {\em local Euler obstruction} coefficients:
\begin{equation}\label{E:Ma-CSM} 
\cMa(X_w) = \sum_{v}e_{w,v}\csm(X_v^{\circ}) \/. 
\end{equation}
Given this, our proof of \Cref{E:main} is based on the following mathematical {\bf Algorithm}:
\begin{enumerate}
\item[(a)] Calculate the Schubert expansion of the CSM classes $\csm(X_w^\circ)$ using the 
recursive algorithm from \S \ref{sec:csm-rec} below, based on
Demazure-Lusztig operators;

\item[(b)] Calculate the Schubert expansion of the Mather classes $\cMa(X_w)$ using the type 
independent formula from \eqref{E:mather-formula} below;

\item[(c)] Calculate the transition matrix $(e_{w,v})$ between the two classes, and compare
to the matrix of $(P_{w,v}(1))$ of Kazhdan-Lusztig coefficients.

\end{enumerate}

We note that the algorithm in (a) holds for any $G/P$, and the one in (b) holds for any cominuscule
$G/P$. 

As a byproduct of our computations we finished the check of the following positivity 
properties, most of which follow from earlier papers in Lie type A \cite{MR1084458},
types A,B,D \cite{boe.fu}, type C \cite{levan.raicu:euler}. This was announced in
\cite{mihalcea.singh:conormal}.

\begin{theorem}\label{thm:positivity} Let $G/P$ be a cominuscule space, let 
$X_w \subset G/P$ be a Schubert variety, and let $v \le w$ in Bruhat order. Then the following positivity properties hold:
\begin{enumerate} \item (Positivity of Euler obstruction) The local Euler obstruction $e_{w,v} \ge 0$ and 
$e_{w,v} > 0$ if $G$ is simply laced. 

\item (Positivity for Mather classes) Consider the Schubert expansion 
\begin{equation*}\cMa(X_w) = \sum_{v \le w} b_{w,v} [X_v] \quad \in H_*(G/P) \/. \end{equation*}
Then $b_{w,v} \ge 0$ and $b_{w,v} > 0$ if $G$ is simply laced.

\item (Strong positivity for CSM classes) Let $G$ be of Lie type $E_6,E_7$ and recall the 
Schubert expansion from \eqref{E:CSM-Schub} above:
\[ \csm(X_w^\circ) = \sum_{v \le w} a_{w,v} [X_v] \quad \in H_*(G/P) \/.\]
Then $a_{w,v} >0$. 

\end{enumerate}
\end{theorem}

The more general strong positivity is conjectured to hold for the CSM class of 
any Schubert cell in any homogeneous space $G/P$. 
Huh's results from \cite{huh:csm} show that each `expected' 
homogeneous component of the CSM classes of Schubert cells 
in Grassmann manifolds is strictly effective, supporting this conjecture. 
The full conjecture seems to be wide open in general.

\section{Proof of \Cref{thm:intro}: theoretical setup}
We prove \Cref{thm:intro} by direct calculation, using the Algorithm from the previous section and a calculation using \texttt{Sage}. We explain next the mathematical foundation of the Algorithm
and in the next section we will explain in more detail the code, and provide its results.

\subsection{Schubert data}\label{sec:schub-data} We follow the notation from \cite{mihalcea.singh:conormal}, and we 
recall it next. Let $G$ be a complex semisimple Lie group with
a maximal torus $T$, and a pair of Borel subgroups $B$ and $ B^-$ satisfying $B\cap B^-=T$.
We denote by  $R$ (resp. $R^+$, $R^-$, $\Delta$) the set of positive (resp. negative, simple) roots.
The set $R$ is equipped with a partial order given by $\alpha < \beta$ if {$\alpha\neq\beta$ and} $\beta - \alpha $ is a non-negative combination of positive roots.
The Weyl group $W:=N_G(T)/T$ associated to $(G,T)$ is a Coxeter group generated by the simple reflections $s_i:=s_{\alpha_i}$, for $\alpha_i \in \Delta$.
Denote by $\ell:W \to \mathbb{N}$ the length function and by $w_0$ the longest element.

Recall that the parabolic subgroups $P\supset B$ are in bijection with the subsets $S \subset \Delta$.
We denote by $R_P^+$ the subset of $R^+$ consisting of roots whose support is contained in $S$.
The Weyl group $W_P$ of $P$ is generated by the simple reflections $s_i$, for $\alpha_i \in S$.
Denote by $w_P$ the longest element in $W_P$, and let $W^P$ be the set of \emph{minimal length representatives} for the cosets in $W/W_P$.
The coset $wW_P$ has a unique minimal length representative $w^P \in W^P$;
we set $\ell(wW_P) := \ell(w^P)$.

The space $G/P$ is a projective manifold of dimension $\ell(w_0W_P)$.
For $w \in W^P$,
the $B$-orbit $X_w^{P,\circ}= B w P/P$, and the $B^-$-orbit $(X^{P})^{w,\circ}= B^- w P/P$,
are opposite Schubert cells for $w$, and there are isomorphisms, 
$X_w^{P,\circ}\simeq \C^{\ell(w)}$ and $(X^P)^{w,\circ}\simeq \C^{\dim G/P - \ell(w)}$.
The \emph{Schubert varieties} $X_w^P$ and $(X^P)^w$ are the closures of the Schubert 
cells $X_w^{P,\circ}$ and $(X^P)^{w,\circ}$ respectively. The fundamental classes 
$\{[X_{w}^P] \}_{w \in W^P}$ form a $\Z$-basis for the homology group $H_*(G/P)$.
Taking the cap product with $[X_{w_0W_P}^P]$ allows us to identify $H_*(G/P)$ with the cohomology
ring $H^*(G/P)$. Under this identification, 
$[X_w^P] \in H_{2 \ell(w)}(G/P) = H^{2(\dim G/P - \ell(w))}(G/P)$.

Every $P$-representation $V$ determines a $G$-equivariant vector bundle, $G \times^P V \to G/P$.
The points of $G\times^PV$ are equivalence classes $[g,v]$,
for pairs $(g,v) \in G \times V$ such that $(g,v) \simeq (gp^{-1},pv)$,
and the $G$-action on $G\times^PV$ is given by left multiplication, $g.[g',v] := [gg',v]$.
The main examples considered in this note are the following.
If $P=B$ is a Borel subgroup, we will take $V:= \C_\lambda$, the one dimensional $B$-module of character $\lambda$.
The resulting line bundle is $\mathcal{L}_\lambda:=G \times^B \C_\lambda$.
Let $\mathfrak{p}$ and $\mathfrak{g}$ be the Lie algebras of $P$ and $G$ respectively, and let $\mathfrak u_P$ be the unipotent radical of $\mathfrak p$.
The spaces $\mathfrak{p}$, $\mathfrak{g}$, and $\mathfrak u_P$ are $P$-stable under the adjoint action.
We have $T(G/P)=G \times^P \mathfrak{g}/\mathfrak{p}$,
and $T^*(G/P)=G\times^P\mathfrak u_P$.

A simple root $\alpha$ is called \emph{cominuscule} if it appears with coefficient $1$ in the highest root of $R$.
Let $\alpha$ be cominuscule, and let $P$ be the parabolic subgroup corresponding to $S=\Delta\backslash\{\alpha\}$.
Then $G/P$ is called a cominuscule space. A 
complete list of cominuscule spaces is given in the introduction. 

\subsection{Calculations of CSM and Mather classes} The key calculation we need to perform
is that of the {\em local Euler obstructions} $e_{w,v}$ from \eqref{E:Ma-CSM} above. 
To calculate the coefficients $e_{w,v}$ we recall algorithms calculating the relevant CSM classes and the Mather classes in a cominuscule $G/P$.  

\subsubsection{CSM classes of Schubert cells}\label{sec:csm-rec} We follow an algorithm proved in \cite{aluffi.mihalcea:eqcsm}, and we refer to {\em loc.~cit.} 
for further references and details. 

Let $w \in W$. If $\lambda$ is a weight, the {\em Chevalley formula} in $H^*(G/B)$ states that 
\begin{equation}\label{E:chevalley}
c_1(\mathcal L_\lambda) \cap [X_w^B] =  \sum \langle -\lambda, \alpha^\vee \rangle [X_{w s_\alpha}^B],
\end{equation}
where the sum is over all positive roots $\alpha$ such that $\ell(w s_\alpha) = \ell(w) -1$.
See e.g. \cite[Thm.~8.1]{buch.m:nbhds} or \cite{brion:flagv}.

Fix a simple root $\alpha_i$ and denote by $P_i$ to be the minimal parabolic group determined by $\alpha_i$. Denote by $p_i:G/B \to G/P_i$ the natural projection. The {\em Bernstein-Gelfand-Gelfand} (BGG) operator
$\partial_i:H^*(G/B) \to H^*(G/B)$ is defined by $\partial_i= p_i^*\circ (p_i)_*$ and it satisfies
\begin{equation}\label{E:BGG} \partial_i [ X_w^B] = \begin{cases} [X_{ws_i}^B] & \textrm{ if } \ell(ws_i) > \ell(w) \/;\\
0 & \textrm{ otherwise} \/. \end{cases} \end{equation}
The relative tangent bundle of $p_i$
is the line bundle $T_{p_i}= \mathcal{L}_{-\alpha_i}$. 
The {\em Demazure-Lusztig} (DL) operator $\mathcal{T}_i: H^*(G/B) \to H^*(G/B)$ is defined by 
\begin{equation}\label{E:DLdef} \mathcal{T}_i = c(T_{p_i}) \partial_i - id = 
(1+ c_1(\mathcal{L}_{-\alpha_i})) \partial_i - id \/. \end{equation} 
The CSM class $\csm(X_w^\circ) \in H_*(G/B)$ may be calculated by using the formula
\begin{equation}\label{E:DL-csm} \mathcal{T}_i(\csm(X_w^{B,\circ})) = \csm(X_{ws_i}^{B,\circ}) \/. \end{equation}
proved in \cite{aluffi.mihalcea:eqcsm}. Together with equations \eqref{E:chevalley} and \eqref{E:BGG},
equation \eqref{E:DL-csm} gives a recursive formula to calculate the Schubert expansion of any CSM class $\csm(X_w^{B,\circ})$,
starting from the initial value
$\csm(X_{id}^B) = [X_{id}^B]$ (the class of the $B$-fixed point in $G/B)$. 

Once the classes in $G/B$ are known, the functoriality of MacPherson's transformation gives
those in $G/P$, as follows. Let $\pi: G/B \to G/P$ be the projection. Then for $w \in W$, 
\begin{equation} \csm(X_{wW_P}^{P,\circ}) = \pi_*(\csm(X_w^{B,\circ})) \quad \in H_*(G/P) \/. \end{equation}
Recall that in terms of Schubert classes $\pi_*[X_w^B]=[X_w^P]$ if $w \in W^P$, and 
$\pi_*[X_w^B]= 0$ otherwise.
\subsubsection{Mather classes of Schubert varieties}\label{sec:mather-form} Let $G/P$ be a cominuscule space. We recall
an algorithm from \cite{mihalcea.singh:conormal} calculating the Schubert 
expansions of Mather classes $\cMa(X_w^P)$. 

For $w \in W^P$, let $I(w)$ denote the inversion set of $w$, 
i.e., the set of positive roots $\alpha$ satisfying $w(\alpha) <0$. Then
it was proved in \cite{mihalcea.singh:conormal} that 
the Mather class of $X_w^P$ is given by 
\begin{equation}\label{E:mather-formula} \cMa(X_w^P) = \pi_* (\prod_{\alpha \in I(w)} c(\mathcal{L}_{-\alpha}) \cap [X_w^B]) \/. \end{equation}
 
\subsection{Irreducibility} We briefly recall 
how the calculation of the local Euler obstructions $e_{w,v}$ gives information about the
multiplicities of the intersection homology sheaves of Schubert varieties. Our approach 
is detailed in \cite{AMSS:shadows} and \cite{mihalcea.singh:conormal}, and goes back to results
in \cite{sabbah:quelques,ginzburg:characteristic}.  

Consider the expansion of the intersection homology sheaf 
\begin{equation}\label{E:ICexp} IH(X_w^P) = \sum m_{w,v} [T^*_{X_v^P} X]  \/.\end{equation} 
This expansion holds in the group of conic Lagrangian cycles in $T^*_{G/P}$, where the adjective conic 
refers to the $\C^*$-dilation action on the cotangent fibres. 
We wish to check that the characteristic cycle giving $IH(X_w^P)$ is irreducible, or, equivalently,
that $IH(X_w^P) = [T^*_{X_v^P} X^P]$. Recall the expansion
\begin{equation}
\label{E:cMa2}
\cMa(X_w^P) = \sum_{v}e_{w,v}\csm(X_v^{P,\circ}) \/.
\end{equation}
It is proved in \cite{boe.fu}, see also \cite{AMSS:shadows}, 
that the characteristic cycle $IH(X_w^P)$ is irreducible if and only if the local Euler obstruction satisfies
\begin{equation}
\label{E:e=P}
e_{w,v}= P_{w,v}(1)
\end{equation} 
for all $v\in W^P$, where $P_{w,v}(q)$ is the Kazhdan-Lusztig polynomial. 
(If the equality does not hold, one can still utilize the Kazhdan-Lusztig conjectures 
to obtain that 
\[ P_{w,v}(1) = \sum_u m_{w,u} e_{u,v} \]
thus recovering the matrix of IH multiplicities from 
the matrices $(e_{w,v})$ and $(P_{w,v}(1))$; 
see \cite{AMSS:shadows}.)


\input{computational-aspects.tex}

\bibliography{conormal}
\bibliographystyle{alpha}

\end{document}

%% file: preamble.tex
\usepackage{amsmath,amsthm,amsfonts,amssymb}
\usepackage{framed}
\usepackage{nicefrac}
\usepackage{color,graphicx}
\usepackage[x11names, rgb]{xcolor}
\usepackage{tikz}
\usetikzlibrary{cd}
\usepackage{hyperref}
\usepackage{pdfsync}
\usepackage{enumitem}
\usepackage{listings}
\usepackage[capitalize,nameinlink]{cleveref}
\crefname{conjecture}{Conjecture}{Conjectures}
\Crefname{conjecture}{Conjecture}{Conjectures}
\usepackage[utf8]{inputenc}
\usepackage{ytableau}
\usepackage[position=bottom]{subfig}

\theoremstyle{plain}
\newtheorem{theorem}{Theorem}

\newtheorem*{theorem*}{Theorem}

\newtheorem{corollary}[theorem]{Corollary}

\crefname{mainTheorem}{Theorem}{Theorems}

\theoremstyle{definition}

\theoremstyle{remark}

\newtheorem*{remark*}{Remark}

\numberwithin{theorem}{section}

%% file: definitions.tex
\def\C{\ensuremath{\mathbb C}}

\newcommand\set[2]{\ensuremath{\left\{#1\mid#2\right\}}}

\newcommand{\csm}{{c_{\text{SM}}}}

\DeclareMathOperator{\Fl}{Fl}

\newcommand{\cMa}{{c_{\text{Ma}}}}

\newcommand{\Z}{{\mathbb Z}}

\DeclareMathOperator{\SL}{{SL}}

\DeclareMathOperator{\Gr}{Gr}

%% file: cotcomin_table-1.tex
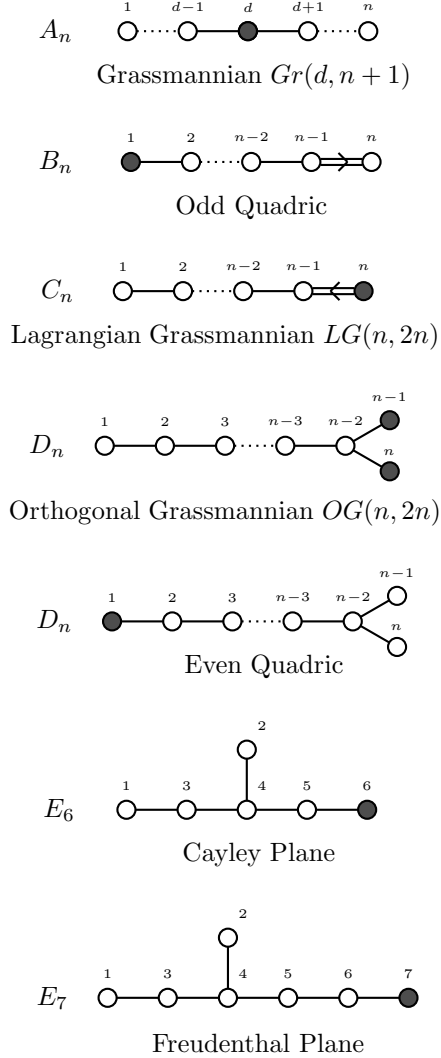
\begin{figure}

\begin{tikzpicture}[scale=.4]
\draw (-1.5,0) node[anchor=east]  {$A_{n}$};
\foreach \x in {0,2,4,6,8}
\draw[thick] (\x cm,0) circle (.3cm);
\draw[thick,fill=black!70] (4 cm,0) circle (.3cm);
\draw[dotted, thick] (0.3 cm,0) -- +(1.4 cm,0);
\foreach \y in {2.3,4.3}
\draw[thick] (\y cm,0) -- +(1.4 cm,0);
\draw[dotted, thick] (6.3 cm,0) -- +(1.4 cm,0);
\draw (0,.8) node {$\scriptscriptstyle{1}$};
\draw (2,.8) node {$\scriptscriptstyle{d-1}$};
\draw (4,.8) node {$\scriptscriptstyle{d}$};
\draw (6,.8) node {$\scriptscriptstyle{d+1}$};
\draw (8,.8) node {$\scriptscriptstyle{n}$};

\draw(4.2, -1.5) node {Grassmannian $Gr(d,n+1)$};
\end{tikzpicture}

\vspace{3mm}

\begin{tikzpicture}[scale=.4]

\draw (13.4,0) node[anchor=east]  {$B_n$};
\draw(19, -1.5) node {Odd Quadric};
\draw (24,0) node[anchor=west]  {\ };
\foreach \x in {17,19,21,23}
\draw[thick,fill=white!70] (\x cm,0) circle (.3cm);
\draw[thick,fill=black!70] (15 cm,0) circle (.3cm);
\draw[dotted,thick] (17.3 cm,0) -- +(1.4 cm,0);
\foreach \y in {15.3,19.3}
\draw[thick] (\y cm,0) -- +(1.4 cm,0);
\draw[thick] (21.3 cm, 1 mm) -- +(14 mm,0);
\draw[thick] (21.3 cm, -1 mm) -- +(14 mm,0);
\draw[thick] (21.9 cm, 3 mm) -- +(3 mm, -3 mm);
\draw[thick] (21.9 cm, -3 mm) -- +(3 mm, 3 mm);
\draw (15,.8) node {$\scriptscriptstyle{1}$};
\draw (17,.8) node {$\scriptscriptstyle{2}$};
\draw (19,.8) node {$\scriptscriptstyle{n-2}$};
\draw (21,.8) node {$\scriptscriptstyle{n-1}$};
\draw (23,.8) node {$\scriptscriptstyle{n}$};
\end{tikzpicture}

\vspace{3mm}

\begin{tikzpicture}[scale=.4]
\draw (-0.8,0) node[anchor=east]  {$C_{n}$};
\draw(3.9, -1.5) node {Lagrangian Grassmannian $LG(n,2n)$};
\foreach \x in {0.5,2.5,4.5,6.5}
\draw[thick,fill=white!70] (\x cm,0) circle (.3cm);
\draw[thick,fill=black!70] (8.5 cm,0) circle (.3cm);
\draw[thick] (0.8 cm,0) -- +(1.4 cm,0);
\draw[dotted,thick] (2.8 cm,0) -- +(1.4 cm,0);
\draw[thick] (4.8 cm,0) -- +(1.4 cm,0);
\draw[thick] (6.8 cm, .1 cm) -- +(1.4 cm,0);
\draw[thick] (6.8 cm, -.1 cm) -- +(1.4 cm,0);
\draw[thick] (7.4 cm, 0 cm) -- +(.3 cm, .3 cm);
\draw[thick] (7.4 cm, 0 cm) -- +(.3 cm, -.3 cm);
\draw (0.5,.8) node {$\scriptscriptstyle{1}$};
\draw (2.5,.8) node {$\scriptscriptstyle{2}$};
\draw (4.5,.8) node {$\scriptscriptstyle{n-2}$};
\draw (6.5,.8) node {$\scriptscriptstyle{n-1}$};
\draw (8.5,.8) node {$\scriptscriptstyle{n}$};
\end{tikzpicture}

\vspace{3mm}

\begin{tikzpicture}[scale=.4]

\draw (14,0) node[anchor=east]  {$D_{n}$};
\draw (19,-2.2) node  {Orthogonal Grassmannian $OG(n,2n)$};
\foreach \x in {17,19,21,23}
\draw[thick,fill=white!70] (\x cm,0) circle (.3cm);
\draw[thick] (15 cm,0) circle (.3cm);
\draw[xshift=23 cm,thick,fill=black!70] (30: 17 mm) circle (.3cm);
\draw[xshift=23 cm,thick,fill=black!70] (-30: 17 mm) circle (.3cm);
\draw[dotted,thick] (19.3 cm,0) -- +(1.4 cm,0);
\foreach \y in {15.3, 17.3, 21.3}
\draw[thick] (\y cm,0) -- +(1.4 cm,0);
\draw[xshift=23 cm,thick] (30: 3 mm) -- (30: 14 mm);
\draw[xshift=23 cm,thick] (-30: 3 mm) -- (-30: 14 mm);
\draw (15,.8) node {$\scriptscriptstyle{1}$};
\draw (17,.8) node {$\scriptscriptstyle{2}$};
\draw (19,.8) node {$\scriptscriptstyle{3}$};
\draw (21,.8) node {$\scriptscriptstyle{n-3}$};
\draw (23,.8) node {$\scriptscriptstyle{n-2}$};
\draw (24.45,1.6) node {$\scriptscriptstyle{n-1}$};
\draw (24.45,-.2) node {$\scriptscriptstyle{n}$};

\end{tikzpicture}

\vspace{3mm}

\begin{tikzpicture}[scale=.4]

\draw (14,0) node[anchor=east]  {$D_{n}$};
\draw (23,-1.5) node[anchor=east] {Even Quadric};
\foreach \x in {17,19,21,23}
\draw[thick,fill=white!70] (\x cm,0) circle (.3cm);
\draw[thick,fill=black!70] (15 cm,0) circle (.3cm);
\draw[xshift=23 cm,thick] (30: 17 mm) circle (.3cm);
\draw[xshift=23 cm,thick] (-30: 17 mm) circle (.3cm);
\draw[dotted,thick] (19.3 cm,0) -- +(1.4 cm,0);
\foreach \y in {15.3, 17.3, 21.3}
\draw[thick] (\y cm,0) -- +(1.4 cm,0);
\draw[xshift=23 cm,thick] (30: 3 mm) -- (30: 14 mm);
\draw[xshift=23 cm,thick] (-30: 3 mm) -- (-30: 14 mm);
\draw (15,.8) node {$\scriptscriptstyle{1}$};
\draw (17,.8) node {$\scriptscriptstyle{2}$};
\draw (19,.8) node {$\scriptscriptstyle{3}$};
\draw (21,.8) node {$\scriptscriptstyle{n-3}$};
\draw (23,.8) node {$\scriptscriptstyle{n-2}$};
\draw (24.45,1.6) node {$\scriptscriptstyle{n-1}$};
\draw (24.45,-.2) node {$\scriptscriptstyle{n}$};

\end{tikzpicture}

\vspace{3mm}

\begin{tikzpicture}[scale=.4]
\draw (-2.4,0) node[anchor=east]  {$E_{6}$};
\draw(3.4, -1.5) node {Cayley Plane};
\draw (8,0) node[anchor=west]  {\ };
\foreach \x in {-1,1,3,5}
\draw[thick,fill=white!70] (\x cm,0) circle (.3cm);
\draw[thick,fill=white!70] (3 cm, 2 cm) circle (.3cm);
\foreach \x in {7}
\draw[thick,fill=black!70] (\x cm,0) circle (.3cm);
\foreach \y in {-0.7, 1.3, 3.3, 5.3}
\draw[thick] (\y cm,0) -- +(1.4 cm,0);
\draw[thick] (3 cm,.3 cm) -- +(0,1.4 cm);
\draw (-1,.8) node {$\scriptscriptstyle{1}$};
\draw (1,.8) node {$\scriptscriptstyle{3}$};
\draw (3.5,.8) node {$\scriptscriptstyle{4}$};
\draw (5,.8) node {$\scriptscriptstyle{5}$};
\draw (7,.8) node {$\scriptscriptstyle{6}$};
\draw (3.5,2.8) node {$\scriptscriptstyle{2}$};
\end{tikzpicture}

\vspace{3mm}

\begin{tikzpicture}[scale=.4]
\draw (12,0) node[anchor=east]  {$E_{7}$};
\draw(18, -1.5) node {Freudenthal Plane};
\foreach \x in {13,15,17,19,21}
\draw[thick,fill=white!70] (\x cm,0) circle (.3cm);
\draw[thick,fill=white!70] (17 cm, 2 cm) circle (.3cm);
\draw[thick,fill=black!70] (23 cm,0) circle (.3cm);
\foreach \y in {13.3, 15.3, 17.3, 19.3, 21.3}
\draw[thick] (\y cm,0) -- +(1.4 cm,0);
\draw[thick] (17 cm,.3 cm) -- +(0,1.4 cm);
\draw (13,.8) node {$\scriptscriptstyle{1}$};
\draw (15,.8) node {$\scriptscriptstyle{3}$};
\draw (17.5,.8) node {$\scriptscriptstyle{4}$};
\draw (19,.8) node {$\scriptscriptstyle{5}$};
\draw (21,.8) node {$\scriptscriptstyle{6}$};
\draw (23,.8) node {$\scriptscriptstyle{7}$};
\draw (17.5,2.8) node {$\scriptscriptstyle{2}$};

\end{tikzpicture}

\vspace{4mm} 
\caption{ Cominuscule spaces and the corresponding Dynkin diagrams with the cominuscule root marked in black; cf.~\cite{bourbaki}.}
\label{TBL:comin}
\end{figure}

%% file: computational-aspects.tex
\section{Proof of Theorem 1.1: Computational Aspects}
Consider formula \eqref{E:mather-formula} computing the Mather class of a Schubert variety $G/P$.
We see that the Mather class is the image of a class $c\in H_*(G/B)$ under the push-forward map $\pi_*:H_*(G/B)\to H_*(G/P)$.
The class $c$ itself is the image of a Schubert class under a sequence of linear operators $c_1(\mathcal L_\alpha)$.
The coefficients of these operators, viewed as matrices, are combinatorial expressions depending on the underlying root system.

In the next sections,
we discuss some of the choices involved in the implementations of \cref{E:chevalley} and \cref{E:DLdef}

\subsection{Choice of Basis}
To compute the Mather (resp. CSM class),
we first to compute the class $\prod c_1(\mathcal L_\alpha) \cap [X_w^B]$ (resp. $\mathcal T_{w^{-1}}[X_w^B]$) in $H_*(G/B)$.

There are at least two natural choices of bases of $H_*(G/B)$ which can be used for this calculation.
One natural choice is the Schubert basis $\set{[X_w^B]}{w\in W}$,
which has the benefit of allowing us to work with integer coefficients.
A second choice of basis is the set of $T$-fixed point classes in the equivariant homology $H_*^T(G/B)$.
Note that the operators $c^T_1(\mathcal L_\alpha)$ act diagonally with respect to this basis,
but with coefficients being rational functions in the simple roots.

The Schubert basis turns out to be the efficient choice here.
The Euler obstructions can be computed from the non-equivariant Mather and CSM classes;
computing the equivariant versions of these classes is slower,
and in this case, unnecessary.
Experimentally, we observed that after a week of computation with unoptimized code,
we were able to compute only $24$ out of the $27$ Mather classes in $E_6$
when working with the point class basis.

\subsection{Representation of Homology Classes and the actions of a Chern classes and DL operators}
\label{homologyCounter}
Given a class $c\in H_*(G/B)$ and the operator $c_1(\mathcal L_\alpha)$,
what is an efficient way to represent the homology class $c$ and compute the action of $c_1(\mathcal L_\alpha)$ and $\mathcal T_i$ on $c$?

A natural idea is to store the class $c$ as a column vector with respect to the Schubert basis,
compute the matrices of the linear operators $c_1(\mathcal L_\alpha)$ and $\mathcal T_i$,
and then evaluate the action of these operators on homology classes via matrix multiplication.
This turns out to be infeasible.
Computing the matrix associated to $c_1(\mathcal L_\alpha)$ is both time-consuming,
and has a prohibitively high memory cost.
A rough empirical estimate suggests that storing all these matrices (even as sparse matrices) will require hundreds of gigabytes of memory in the $E_7$ case.

Instead, we store the homology class $c$ as a Counter object (see \cite{Counter,Dict}),
with Weyl group elements being the keys and the corresponding coefficients being the values.
We call a function \texttt{lUpdate} to update the value of a homology class under the action of a Chern class operator $c_1(\mathcal L_\alpha$)
and a function \texttt{tOp} to compute the action of $\mathcal T_i$ on the homology class stored in a Counter.
The trade-off in this choice is that while the entries of the matrix of $c_1(\mathcal L_\alpha)$ and $\mathcal T_i$ are calculated each time these functions are called,
only the entries corresponding to non-zero coefficients in the class $c$ are calculated.

Observe that for $R'$ any subset of $R_P^-$,
the class $\left(\prod\limits_{\alpha\in R'} c_1(\mathcal L_\alpha)\right)\cap [X_w^B]$
is supported on a subset of the set $\set{[X_v^B]}{v\leq w}$.
Let $w_0^P$ denote the minimal representative in $W^P$ of the longest element $w_0$.
A computer verification shows that 
the cardinality of the set $\set{[X_v^B]}{v\leq w^P_0}$ is significantly smaller than $|W|$,
see Table 1.
This explains the better performance of our method;
even though the calculation of some coefficients is repeated,
the calculation of many other coefficients is avoided all-together,
since we only compute the action of the operators $c_1(\mathcal L_\alpha)$ and $\mathcal T_i$ on the classes $[X_v^B]$,
for $v\leq w_0^P$.

\begin{table}[ht]
\begin{tabular}{|c|c|c|c|c|}
\hline
\label{sizeTable}
Dynkin diagram & Cominuscule Space & $|W|$ & $|W^P|$& $\#\set{v\in W}{v\leq w_0^P}$\\ \hline
$E_6$ & Cayley Plane& $51840$ & $27$ & $5264$\\ \hline
$E_7$ & Freudenthal Plane &$2903040$ & $56$ & $228696$ \\ \hline
\end{tabular}
\end{table}

\subsection{Computing the Kazhdan-Lusztig Polynomial}
Fokko du Cloux's Coxeter3 program is considered state of the art for the computation of Kazhdan-Lusztig polynomials.
As such, it was the natural choice for generating the Kazhdan-Lusztig classes in the cominuscule Grassmannian.
However, the Coxeter3 library in SageMath is not optimized for hashing.
In particular, storing a homology class as a Counter object (see \cref{homologyCounter}) is not feasible with Coxeter3,
because the space requirements of Coxeter3 objects are huge.
Therefore, we use the default implementation of Coxeter groups in Sage for the computation of the Mather and CSM classes,
and the Coxeter3 implementation for the computation of the Kazhdan-Lusztig classes.

\subsection{The Code}
We present key parts of the code with some explanations.
\texttt{ 
\begin{lstlisting}[language=python]
54  W = CoxeterGroup( lie_type + str(rank))
55  S = [W.simple_root_index(i) for i in W.index_set()]
56  alpha = [W.roots()[i] for i in S]
57  alpha = {i+1: root for i, root in enumerate(alpha)}
58  S_P = list(range( 1, d)) + list(range( d + 1, rank + 1))
59  W3 = CoxeterGroup(lie_type + str(rank),
                      implementation='coxeter3')
\end{lstlisting} }
The variable \texttt{W} stores an instance of the Weyl group using the standard implementation of a Coxeter group in Sage,
while \texttt{W3} stores an instance of the Weyl group in the Coxeter3 implementation.
The variable \texttt{alpha} is a dictionary mapping the integers $1,\cdots,n$ to the corresponding simple root.
The variable \texttt{S\_P} is the list of simple roots contained in $P$.
\texttt{ 
\begin{lstlisting}[language=python]
61  def cartan_matrix(W):
64      C = [[v for v in row] for row in W.coxeter_matrix()]
65      for i in range(len(C)):
66          for j in range(len(C[0])):
67              x = C[i][j]
68              C[i][j] = 2 if x == 1 else 0 if x == 2 else -1
69      return Matrix(C)

71  C = cartan_matrix(W)

82  def pairing(a,b): 
86      return a*C*b
\end{lstlisting} }
The Chevalley formula invokes the bilinear pairing $\left\langle\,,\,\right\rangle$,
which we can compute using the Cartan matrix.
Since we only deal with simply laced Dynkin diagrams,
the Cartan matrix can be deduced from the Coxeter matrix of $W$.
The function \texttt{cartan\_matrix} computes the Cartan matrix,
and the function \texttt{pairing} uses \texttt{cartan\_matrix} to compute the bilinear pairing $(\alpha,\beta)\mapsto (\alpha^\vee,\beta)$.
\texttt{ 
\begin{lstlisting}[language=python]
76  root_of_reflection = {
77      W.reflections()[root]: root\
78          for root in W.roots()\
79          if root > 0
80  }
\end{lstlisting} }
The Weyl group \texttt W has an inbuilt method to compute the reflection across a given root.
Since no method is provided to compute the reverse mapping directly,
we store this reverse mapping in the variable \texttt{root\_of\_reflection}.
\texttt{ 
\begin{lstlisting}[language=python]
88  def generateWP():
91      new = {W.one()}
92      by_len = []
93      while new:
94          by_len.append(new)
95          new = set()
96          for w in by_len[-1]:
97              for sref in W.simple_reflections():
98                  v = sref * w
99                  if v.length() > w.length() and\
                        all( v * alpha[i] > 0 for i in S_P):
100                     new.add(v)
101     new = []
102     for sub in by_len:
103         for w in sub:
104             new.append(w)
105     return new

107 WP = generateWP()
\end{lstlisting} }
The function \texttt{generateWP} returns a list containing the elements of $W^P$.
A naive way to generate $W^P$ involves looping over the elements of $W$
and filtering in the elements which are minimal with respect to $W_P$,
for example by checking if $w(\alpha)>0$ for all $\alpha\in S_P$.

While this method is computationally feasible, there is a faster way.
Recall that if $s_{\alpha_1}\cdots s_{\alpha_i}$ is minimal with respect to $W_P$,
so is $s_{\alpha_2}\cdots s_{\alpha_i}$.
This allows us to build the set $W^P$ inductively.
We start with the identity element.
At each step, having computed all the elements of $W^P$ of a particular length,
we generate all elements in $W$ of length one more,
and filter in those which are minimal with respect to $W_P$.

This method has the further incidental benefit that the generated list is sorted in increasing order by Coxeter length.
In particular, the Mather and CSM classes will be upper-triangular
with respect to the ordered basis corresponding to the list \texttt{WP = generateWP()}.
\texttt{ 
\begin{lstlisting}[language=python]
110 def lUpdate( root, vec):
117     c = deepcopy(vec)
118     for w in c:
119         for v in w.lower_cover_reflections():
120             p = -pairing(root, root_of_reflection[v])
121             vec[w*v] += p * c[w]
\end{lstlisting} }
The function \texttt{lUpdate} takes as input a root $\alpha$ and a homology class $c$ (the variable \texttt{vec}).
It modifies \texttt{vec} to the homology class $c_1(\mathcal L_\alpha) \cap c$.
\texttt{ 
\begin{lstlisting}[language=python]
123 def tOp( k, vec):
129     res = Counter()
130     for w, cof in vec.items():
131         if w * alpha[k] < 0: 
132             res[w] -= cof
133             continue
134         res[w] -= cof
135         v = w * W.simple_reflection(k)
136         res[v] += cof
137         for u in v.lower_cover_reflections(): 
138             beta = root_of_reflection[u]
139             res[v*u] += cof * pairing( beta, alpha[k])
140     return res
\end{lstlisting} }
The function \texttt{tOp} takes as input an index $k$ and a homology class $c$ (the variable \texttt{vec}),
and returns the homology class $\mathcal T_{s_k}\cdot c$, leaving the Counter \texttt{vec} unchanged.
\texttt{ 
\begin{lstlisting}[language=python]
148 def compute_mather_classes():
152     Rn = [root for root in W.roots()\
             if root < 0 and d-1 in root.support()]
156     cMa=[]
157     for w in WP:
158         v = Counter([w])
159         for beta in Rn:
160             if w * beta > 0: 
161                 lUpdate(beta, v)
162         cMa.append([])
163         for u in WP:
164             cMa[-1].append(v[u])
165     return cMa
\end{lstlisting} }
In the code snippet above,
lines 158--161 compute the class $\prod\limits_{\alpha\in R_P^+}c_1(L_\alpha)\cap [X_w^B]$ (in $H^*(G/B)$),
and lines 162--164 append its projection to the list \texttt{cMa}.
\texttt{ 
\begin{lstlisting}[language=python]
167 def compute_csm_classes():
172     csm=[]
173     for w in WP:
174         v = Counter([W.one()])
175         for k in w.reduced_word():
176             v = tOp(k,v)
177         csm.append([])
178         for u in WP:
179             csm[-1].append(v[u])
180             assert csm[-1][-1] >= 0
181             assert u.bruhat_le(w) == (csm[-1][-1] > 0) 
182     return csm
\end{lstlisting} }
In the code snippet above, the loop in lines 175--176 computes $\mathcal T_{w^{-1}}\cdot[X_{id}^B]$ (in $H^*(G/B)$),
and the loop in lines 178--181 appends the projection of this class (to $H^*(G/P)$) to the list \texttt{csm}.
Lines 180--181 verify the strong positivity conjecture.
\texttt{ 
\begin{lstlisting}[language=python]
210 def compute_kl():
211     W3P = [W3.from_reduced_word(x.reduced_word())\
                for x in W_P]
212     w_P = W3.from_reduced_word(longest_element_in_WP())
213     WPmax = [x*w_P for x in W3P]
214     klp = []
215     for u in WPmax:
216         klp.append([])
217         for v in WPmax:
218             poly = W3.kazhdan_lusztig_polynomial( u, v)
219             klp[-1].append( sum(poly.coefficients()))
220     return Matrix(klp)
\end{lstlisting} }
The variable \texttt{W3P} is a list containing the elements of $W^P$ as Coxeter3 objects,
in the same order as the list \texttt{WP}.
The variable \texttt{WPmax} is a list of the maximal representatives.
\texttt{ 
\begin{lstlisting}[language=python]
223 cMaMatrix = Matrix(compute_mather_classes()).transpose()
226 csmMatrix = Matrix(compute_csm_classes()).transpose()
228 eulerObsMatrix = csmMatrix.inverse() * cMaMatrix
231 klp = compute_kl()
\end{lstlisting} }
Lines 223 and 226 compute matrices whose columns are the Mather and CSM classes of Schubert varieties in $G/P$
in the basis \texttt{WP}.
Using \cref{E:cMa2}, we compute the matrix of Euler obstructions,
which we compare with the matrix of Kazhdan-Lusztig polynomials evaluated at $1$.